\documentclass[a4paper,12pt]{article}
\usepackage{amssymb,amsthm,amsfonts}
\newtheorem{theorem}{Theorem}

\newcommand{\Z}{{\mathbb Z}}
\newcommand{\Zpi}{{\mathbb Z}_{p^\infty}}
\newcommand{\fb}{Fern\'andez-Bret\'on}

\newcommand{\Q}{{\mathbb Q}}

\newcommand{\erdos}{{Erd\H os}}

\title{Monochromatic Sumsets in Countable Colourings of Abelian Groups}

\author{Imre Leader \and Kada Williams}
\newcommand{\Addresses}{{
    \bigskip\footnotesize

    Imre~Leader, Department of Pure Mathematics and Mathematical Statistics, Centre for Mathematical
    Sciences, Wilberforce Road, Cambridge CB3 0WB, UK. \\ \nopagebreak
    \textit{Email}: {i.leader@dpmms.cam.ac.uk}

    \medskip

    Kada Williams, Department of Pure Mathematics and Mathematical Statistics, Centre for
    Mathematical Sciences, Wilberforce Road, Cambridge CB3 0WB, UK. \\ \nopagebreak
    \textit{Email}: {kkw25@cam.ac.uk}

    }}

\begin{document}
\maketitle

\begin{abstract}
Fern\'andez-Bret\'on, Sarmiento and Vera showed that whenever a direct sum of sufficiently many
copies of $\Z_4$, the cyclic group of order 4, is countably coloured there are
arbitrarily large finite sets $X$ whose sumsets $X+X$ are 
monochromatic. They asked if the elements of order 4 are necessary, in the
following strong sense: if $G$ is an abelian group having no elements of
order 4, is it always the case there there is a countable colouring of $G$
for which there is not even a monochromatic sumset $X+X$ with $X$ of size 2?
Our aim in this short note is to show that this is indeed the case.
\end{abstract}

\begin{section}{Introduction}
It follows directly from Ramsey's theorem that, whenever the 
naturals are finitely coloured, there is an infinite set $X$ such that all 
pairwise sums of distinct elements of $X$ have the same colour. However, if one asks 
for a stronger conclusion, that the entire sumset $X+X = \{x+y: x,y \in X\}$ (in other
words, including terms of the form $x+x$) is 
monochromatic, then it is easy to see that the answer is no -- indeed, it 
is a simple matter to find a 3-colouring yielding no such monochromatic set.

It is worth mentioning that, surprisingly, it is unknown as to whether or not
this can be achieved with a 2-colouring: this is called Owings' problem
\cite{owings}. For background on this, see \cite{hls}.

If one extends from the naturals to the rationals, the answer remains no: 
there is a 12-colouring of $\Q$ with 
no infinite sumset monochromatic. What about for the reals?
Hindman, Leader and Strauss \cite{hls} 
showed that, for every rational vector 
space of dimension smaller that $\aleph_\omega$, there is a finite 
colouring without an infinite monochromatic sumset. This 
establishes the answer for the reals if we assume CH. In the other
direction, Leader and Russell \cite{lr} showed that there \textit{are} vector spaces where
the answer is yes: if we finitely colour a rational vector space of dimension
$\beth_\omega=\sup\{\aleph_0,2^{\aleph_0},2^{2^{\aleph_0}},\ldots\,\}$ then there
is always an infinite monochromatic sumset. Subsequently, Komj\'ath, Leader, Russell,
Shelah, Soukup and Vidny\'anszky \cite{six} proved that this
positive result also holds for the reals, under a `large cardinal' assumption, and finally Zhang
\cite{zhang} completed the picture by removing this assumption: he gave an ingenious argument to
show that it is consistent that whenever the reals are finitely coloured there is an infinite
monochromatic sumset.

Several of these results have been considered and generalised to finite colourings of abelian groups,
not just rational vector spaces. See for example \fb~\cite{f} and \fb~and Lee \cite{fl}. For general
background and several new results in this direction see \fb, Sarmiento and Vera \cite{fb}.

Turning to colourings with countably many colours, there are results about monochromatic structures when
we countably colour the finite subsets of a (large) set, going back at least as far as
Komj\'ath and Shelah (\cite{ks} and \cite{ks2}) -- see also Komj\'ath \cite{ko} and
\fb~and Rinot \cite{fr} for results about abelian groups. However, for sumsets, the situation is very
different to the descriptions above. 
Indeed, one can never guarantee an infinite sumset, however large the abelian group. This was
proved in \cite{fb}. But what happens for finite sumsets? The authors of \cite{fb} gave an
elegant argument, based on the \erdos-Rado theorem \cite{er}, to
show the following: if we countably colour the direct sum of sufficiently many copies
of $\Z_4$, the cyclic group of order 4, then for any $n$ there is a set $X$ of size $n$ whose
sumset $X+X$ is monochromatic. Here `sufficiently many' is in fact $\beth_\omega$ copies.

They asked if these elements of order 4 are necessary, in the following very strong sense: if $G$
is an abelian group with no elements of order 4, does there always exist a countable colouring of
$G$ for which there does not even exist a set $X$ of size 2 with $X+X$ monochromatic? In other
words, there do not exist distinct $x$ and $y$ such that the elements $2x,2y,x+y$ all have the same
colour. If true, this would make a very attractive counterpart to the result about direct sums
of $\Z_4$. They proved this (in the affirmative) in several cases, such as when $G$ is a torsion
group (meaning that it does not have any elements of infinite order).

Our aim in this note is to prove this in general. 

\begin{theorem}\label{main} 
  Let $G$ be an abelian group not containing any elements of order 4. Then there is a countable
  colouring of $G$ for which there do not exist distinct $x$ and $y$ with $\{ 2x,2y,x+y \}$
  monochromatic.
\end{theorem}

To end this section, we give a brief explanation of the main ideas of the proof. A standard
way to deal with general abelian
groups (see for example \cite{kap}) is to embed them into direct sums of the groups $\Q$ and $\Zpi$
for any
prime $p$ ($\Zpi$ being the Pr\"ufer group, consisting of all $p^k$-th complex roots of unity, for
all $k$).
This is achieved by first embedding our group into a divisible group (meaning that every group element
$g$ is divisible, i.e.~may be divided by $n$ for every positive integer $n$) and then
showing that any divisible group may
be expressed as such a direct sum.

Now, if our starting group $G$ (with no elements of order 4) has no elements of infinite order, then
it turns out that one can
embed $G$ into a `sufficiently nice' group, meaning a group in which all elements not of order 2
are divisible, without creating any elements of order 4. However, if we have elements
of infinite order then this cannot work. Indeed, if $x$ and $y$ have infinite order and their
difference has order 2, and
we assume that we have extended our group by adding elements $g$ and $h$ such that $2g=x$ and
$2h=y$, then the element $g-h$ has order 4.

So instead we have to work with the infinite-order elements separately. We first embed the
odd-order torsion part of the group into a divisible group $D$, and this we will be able to
handle (in the sense of giving some countable colourings) by the structure theorem above.
The group $G$ now lives inside $D \oplus X$ for some $X$, because by Baer's theorem \cite{baer}
a divisible subgroup is always a direct summand, which means that we may `deal with these two
parts separately'. Now, the torsion part $T$ of $X$ is a vector space over
$\Z_2$. But, because we are forced to work with $T$ and not some
divisible extension of $T$, we do not have $T$ as a direct summand of $X$. So instead we look
at the quotient $X/T$, which is torsion-free and so embeds into a direct sum of copies
of $\Q$. This we will again be able to handle. This leaves us with $T$, but unfortunately embedded
merely as a
subgroup, not a direct summand. However, it turns out that a very convenient property of the structure
will allow us to complete the colouring.

Our notation is standard. See \cite{kap} or \cite{fuchs} for general background on abelian groups, including
divisible groups and Baer's theorem.

\end{section}

\begin{section}{Proof of Main Result}

Whenever we have an element of a direct sum $\bigoplus_{i \in I} G_i$, we will always assume that
the index set $I$ is linearly ordered. This means that each
element $x$ of this direct sum has a {\it support}, meaning
$\{ i \in I : x_i \neq 0 \}$, and also a {\it profile}, meaning the ordered finite sequence of the
non-zero $x_i$. Note that if each $G_i$ is countable, and there are only countably many different
groups that appear as a $G_i$ (whatever the cardinality of $I$), then there are only countably
many profiles.

\begin{proof}[Proof of Theorem~\ref{main}]
Given our abelian group $G$ with no elements of order 4, we will build up our colouring in
stages, at each time obtaining stronger and stronger information about a supposed monochromatic
$\{ 2a,2b,a+b \}$, until we obtain a contradiction. Formally we will be taking one colouring, then
taking the product colouring with the next colouring, and so on.

Let $x$ be any element of $G$ of odd order.
We can already divide $x$ by 2 (because $x$ has odd order), but if there is an odd
prime $p$ for which we cannot divide $x$ by
$p$ then we extend $G$ by adding an element that is $x$ divided by $p$. This is a standard
operation: we form $G \times \Z$ and quotient by the subgroup generated by $(x,p)$. The 
resulting group contains a copy of $G$ (namely $G \times \{ 0 \}$) in which we can divide $x$ by $p$, 
because $p(0,-1)=(x,0)$. Note also that no elements of order 4 have been introduced: indeed, if
$4(g,r)=0$ then $4r$ must be a multiple of $p$, whence $r$ is a multiple of $p$, so that
$(g,r)$ was already an element of the copy of $G$, a contradiction.

It is now a standard exercise to `repeat' the above operation (one can perform an iteration,
or a transfinite induction)
to obtain a group $H$, containing $G$ as a subgroup, in which every element of
odd order is divisible. It is of course now sufficient to countably colour $H$ in such a way
that there is no monochromatic $\{2a,2b,a+b \}$.

Let $D$ be the subgroup formed by all the odd-order elements of $H$. This is a
divisible subgroup, by the above, and so by Baer's theorem we can write $H$ as $D \oplus X$ for some
subgroup $X$. Note also
that, because $D$ is divisible, it is a direct sum of copies of the groups $\Zpi$ for odd
primes $p$ -- this is the standard argument mentioned in the introduction (see \cite{kap}).

We start by colouring $H$ by giving $(d,x)$ a colour that is the profile of $d$ -- recall that
$D$ has only countably many profiles. We claim that if $a=(d,x)$ and $b=(e,y)$ have the
property that $\{ 2a,2b,a+b \}$ is monochromatic then in fact $d=e$. To see this, note that we
know that $2d,2e,d+e$ all have the same profile. Now, if the supports of $d$ and $e$ do not
end in the same place then this is impossible, as if the final entry of the profile of $d$ is $t$
and $d$ ends later than $e$ then the profile of $2d$ ends with $2t$ while the profile of $d+e$
ends with $t$. (Note that for this it is vital that $D$ consists only of odd-order elements, so that
$2t$ cannot be zero.) So the supports of $d$ and $e$ end in the same place, and now we may look at
the second-last elements of their supports. These must be equal by the same argument.
Continuing, we see that $d$ and $e$ have the same support. Hence $d=e$: since $2d$ and
$2e$ have the same profile it follows that $d$ and $e$ have the same profile.

Now let $T$ be the torsion subgroup of $X$. Every element of $T$ has order 2, since other orders
would yield either an element of $X$ of order 4 or an element of $X$ of odd order, both of which
are impossible. Moreover, the group $X/T$ is torsion-free, and so embeds into a group $Y$ that is
a direct sum of copies of $\Q$. We now additionally colour by giving $(d,x)$ a colour which
is the profile of $\pi(x)$, where $\pi$ denotes the projection from $X$ to $X/T$, with the
latter viewed as a subgroup of $Y$. By exactly the same argument as above, if
$a=(d,x)$ and $b=(e,y)$ have $\{ 2a,2b,a+b \}$ monochromatic then we must have $\pi(x)=\pi(y)$ (the
same argument works because, as with $D$, $Y$ has no elements of even order). 

To summarise, we now know that $d=e$ and that $x$ and $y$ belong to the same coset of $T$ in $X$,
in other words that $x+T=y+T$. It is here that we are rather fortunate. We claim that in each coset of
$T$ there is at most one element that can be halved. Indeed, suppose that $u+T=v+T$ and that there
exist $g$ and $h$ with $2g=u$ and $2h=v$. Then $2(g-h)$ is a non-zero element of $T$, and so has order
2, which implies that $g-h$ has order 4.

So let us add a 2-colouring of $H$, by the answer to the question `can this element be halved?'.
Now, $2a$ can certainly
be halved, but $2x$ and $x+y$ are distinct elements of the same coset of $T$ (distinct because
otherwise we would have $d=e$ and $x=y$, whence $a=b$). So, by the above
observation, $x+y$ cannot be halved. Thus $a+b$ is not the same colour as $2a$, and we are done.

\end{proof}

It would be very interesting to know what dimension is needed in the result from \cite{fb} about
direct sums of copies of $\Z_4$. Perhaps $\beth_\omega$ is sharp, if we assume CH? It would also
be nice to find out just how many elements of order 4 need to be present to ensure a sumset
$X+X$ with $X$ of size 2 -- in other words, to determine precisely which abelian groups have the property
that whenever they are countably coloured there is such a monochromatic sumset.

\end{section}

\Addresses

\end{document}